\documentclass[reqno]{amsart}

\usepackage{amsmath}
\usepackage{amssymb, latexsym}

\newcommand{\HW}{\textrm{Conv}(W{\mu})}
\newcommand{\no}{{\bar \nu}_1}
\newcommand{\nt}{{\bar \nu}_2}
\newcommand{\nr}{{\bar \nu}_r}
\newcommand{\nik}{{\bar \nu}_k}
\newcommand{\nM}{\nu_M}

\newcommand{\nk}{{\bar \nu}_k}
\newcommand{\nkk}{{\bar \nu}_{k+1}}

\newcommand{\Nk}{\sigma(k)}

\newtheorem{lem}{Lemma}[section]
\newtheorem{thm}[lem]{Theorem}
\newtheorem*{KR}{Proposition [\cite{KR02}, Prop. 4.6]}
\newtheorem*{MT}{Theorem}
\newtheorem*{CM}{Converse to Mazur's Inequality}
\newtheorem{fact}{Fact}

\numberwithin{equation}{section}

\begin{document}

\title{A Converse to Mazur's inequality \\
for split classical groups}
\author{Catherine Leigh}
\address{Department of Mathematics\\
         University of Chicago\\ 
         5734 South University Avenue\\
         Chicago, IL 60637} 
\email{cleigh@math.uchicago.edu}
\date{November 19, 2002}
\maketitle
\section*{}
Our goal is to prove a converse to Mazur's inequality for split classical 
groups.  This work stems from results in~\cite{KR02}, whose notation we will 
follow.  Let $G$ be a split connected reductive group with Borel subgroup $B$
 and let $T$ be  a maximal torus contained in $B$.  We abbreviate $X_* T$ 
to $X$.  We denote by $\mathfrak{a}$  the real vector space 
$X \otimes_{\mathbb Z} {\mathbb R}$, and by $\mathfrak{a}_{\textrm{dom}}$ 
the cone of dominant elements in $\mathfrak{a}$.
For $x,y \in~\mathfrak{a}$ we say $x\leq y$ if $\langle x,\omega \rangle \leq
\langle y,\omega \rangle$ for all fundamental weights $\omega$ and $y-x$ is 
in the linear span of the coroots of $G$.
We write $X_G$ for the quotient of $X$ by  the coroot lattice for $G$, and
$\varphi_G : X \rightarrow X_G$ for the natural projection map.

Now let $\mu \in X$ be $G$-dominant.  The Weyl group $W$ acts on $X$. 
We define the subset
$P_{\mu} \subset X$ by

\[ P_{\mu}  := \left\{ \nu\in X : \begin{array}{c@{\ }l} (i) & 
\varphi_G(\nu)= \varphi_G(\mu),\textrm{ and} \\
     (ii) & \nu \in \HW \end{array} \right\} \]
where $W{\mu}= \{w(\mu): w \in W  \}$ and $\HW$ denotes
the convex hull of $W{\mu}$  in $\mathfrak{a}$.

Let $P$ be a parabolic subgroup of $G$ which contains $B$ and let $M$ 
be the unique Levi subgroup of $P$
containing $T$.
Thus $T$ is also a maximal
torus in
$M$.  We write $X_M$ for the quotient of $X$ by the coroot lattice for $M$,
and let $\varphi_M : X \rightarrow X_M$ be the natural projection map.
Since the coroot lattice for $M$ is a subgroup of the coroot lattice
for $G$, the map $\varphi_G$ factors through $X_M$ via $\varphi_M$.

Now let $F$ be a finite extension of ${\mathbb Q}_p$ and let $L$ be 
the completion
of the maximal unramified extension of $F$.  Denote the ring of integers in
$F$, resp.~$L$, by $\mathcal{O}_F$, resp.~$\mathcal{O}_L$, and set 
$\tilde{K} = G(\mathcal{O}_L)$.
Let $\pi\in \mathcal{O}_F$ be a uniformizer and let $\sigma$ be the relative
Frobenius automorphism of $L/F$.
Recall that $b_1, b_2 \in G(L)$ are said to be $\sigma$-conjugate if 
there exists $g\in G(L)$
such that $b_2 = g^{-1} b_1 \sigma(g)$  and that $B(G)$ denotes the 
set of $\sigma$-conjugacy classes in $G(L)$.
Finally, as in Section 6 of~\cite{Kot97}, define
$B(G, \mu) := \{ b\in B(G): \kappa_G(b)=\varphi_G(\mu) \textrm{ and } 
\bar\nu(b)\leq\mu \}$,
where $\bar\nu(b)\in \mathfrak{a}_{\textrm{dom}}$ is the image of $b$ 
under the Newton map (cf.~\cite{RR96},~\cite{Kot97}) and
$\kappa_G$ is the map (4.9.1)  of~\cite{Kot97}.
We can now state our main theorem.

\begin{CM}
Let $G$ be a split connected reductive group over $F$ which is weakly 
classical, in the sense that each irreducible component of its Dynkin
diagram is  of type $A_n$, $B_n$, $C_n$, or $D_n$.   Let $\mu\in X$ be 
$G$-dominant and let  $b\in B(G, \mu)$.
Then the $\sigma$-conjugacy class of $b$ in $G(L)$ meets
$\tilde{K} \pi^{\mu} \tilde{K}$.
\end{CM}

The converse to Mazur's inequality is proven in~\cite{KR02} 
(resp.~\cite{FR02}), for $G=GL_n$ and $G=GSp_{2n}$ (resp.~$G=GL_n$).
The reader who would like to know how the theorem relates to the reduction
 modulo $p$ of Shimura varieties should consult the survey 
article~\cite{Rap02}.

To prove the converse for all split classical groups, we use the following 
result of~\cite{KR02}.

\begin{KR}
Let $b \in M(L)$ be basic, and let $\mu \in X$.
The $\sigma$-conjugacy class of $b$ in $G(L)$ meets
$\tilde{K}\pi^{\mu}\tilde{K}$ if and only if
$\kappa_M(b) \in \varphi_M(P_{\mu})$.
\end{KR}

We will show that $\varphi_M(P_{\mu})$ is equal to a certain subset 
of $X_M$, described below.  This is the method used in~\cite{KR02} 
for $G=GL_n$ and $G=GSp_{2n}$.
We will provide proofs  for $G=SO_{2n+1}$ and $G=SO_{2n}/\{ \pm 1 \}$ as 
well as an alternate proof for $G=GL_n$.  While we have not yet
proven the result for exceptional groups, we suspect that it is true.

Now we describe the subset of $X_M$ which we will prove is equal to
$\varphi_M(P_{\mu})$.  Recall that 
$\mathfrak{a} = X \otimes_{\mathbb Z} {\mathbb R}$
and define  $\mathfrak{a}_M = X_M \otimes_{\mathbb Z} {\mathbb R}$.  There is
 a natural projection map $pr_M : \mathfrak{a} \rightarrow 
\mathfrak{a}_M$ induced by $\varphi_M$.
Note that $W$ acts on $\mathfrak{a}$ since it acts on $X$.  Therefore
the Weyl group of $M$, denoted by $W_M$,  acts on $\mathfrak{a}$.
The restriction of $pr_M$ to $\mathfrak{a}^{W_M}$ identifies
$\mathfrak{a}^{W_M}$ with $\mathfrak{a}_M$.  With this identification, 
it is easy to see that
$\HW\cap \mathfrak{a}^{W_M} =pr_M(\HW)$.
We can now formulate the theorem.

\begin{MT}
Let $G$ be a split connected reductive group over $F$ which is weakly 
classical in the
above sense and let $M$ and $\mu$ be as above.
Then \[
\varphi_M(P_{\mu}) = \left\{ \nu_1 \in X_M : \begin{array}{c@{\ }l}
(i) & \nu_1 , \mu  \textrm{\ have the same image in }X_G, \textrm{and} \\
     (ii) &  \textrm{the image of } \nu_1 \textrm{ in } \mathfrak{a}_M 
\textrm{\ lies in }
\textup{Conv}(W\mu)\cap \mathfrak{a}^{W_M}
\end{array} \right\}. \]
\end{MT}
For all $G$ we can easily check that the left
hand side is a subset of the right hand side
by  verifying that for $\nu\in P_{\mu}$, the image  $\varphi_M(\nu)$ 
satisfies
conditions $(i)$ and $(ii)$.
First, $\nu$ and $\mu$ have the same image
in $X_G$, hence $\varphi_M(\nu)$ and $\mu$ also have the same image in $X_G$.
Second,  the image of $\varphi_M(\nu)$ in $\mathfrak{a}_M$ is $pr_M(\nu)$,
  which lies in $pr_M(\HW)$
  since $\nu\in \HW$.

To prove that the right hand side is a subset of the left hand side, we
must show that given $\nu_1 \in X_M$ satisfying conditions $(i)$ and $(ii)$,
we can find $\nu\in X$ such that
  $\varphi_M(\nu) = \nu_1$ and $\nu \in \HW$.
Note that if  $\varphi_M(\nu) = \nu_1$, then $\nu$
and $\mu$ will have the same image in $X_G$
by condition $(i)$.
To show that we can find such a $\nu$ for classical
groups, we will examine each type separately.  Before doing so, we will make
several reductions.

We claim that we may assume without loss of generality that $\nu_1$ 
is $G$-dominant.
Indeed, in the theorem only the Weyl group orbit of $\mu$ plays a 
role, and therefore we are free to choose
the Borel subgroup $B$ such that $\nu_1$ is $G$-dominant.
We now do so, and then (cf.~e.g.~Lemma 2.2 (ii) in~\cite{RR96}) the 
condition that $\nu_1\in\HW$ is equivalent to $\nu_1\leq\mu$.

A coweight  $\nu \in X$ is said to be  $G$-minuscule if for
every root $\alpha$ of $G$, we have 
$\langle\alpha,\nu\rangle\in\{ -1, 0, 1 \}$.
We define $M$-minuscule analogously.  Bourbaki (cf.~Ch.~VIII, \S 7 
Prop.~8) shows
that given $\nu_1\in X_M$ there exists a unique $M$-dominant, $M$-minuscule
coweight $\nu\in X$ such that $\varphi_M(\nu)=\nu_1$.  We will show that this
choice of $\nu$ satisfies the desired condition, i.e. that $\nu\in\HW$.

Denote $pr_M(\nu)$ by $\nM$.  We are reduced to proving 
the following:
Given  $\nu, \mu \in X$ such that $\mu$ is $G$-dominant, $\nu$ is
$M$-dominant and $M$-minuscule, and $\varphi_G(\nu)=\varphi_G(\mu)$,
and given that $\nM \in \mathfrak{a}^{W_M}$ is $G$-dominant and $\nM\leq\mu$,
it follows that $\nu\in \HW$.

As we have already stated, Kottwitz and Rapoport proved this theorem for
$G=GL_n$ and $G=GSp_{2n}$.  We will prove the theorem for $G=SO_{2n+1}$
and $G=SO_{2n}/ \{ \pm 1 \}$ as well as giving an alternate proof for $GL_n$.
We claim that to prove the theorem for all split connected reductive 
groups which are weakly classical,
it is enough to do so for the above groups $G$.
We will need four easy facts:
\begin{fact}
If the theorem holds for $G$ and the center of $G$ is connected, then 
the theorem holds for the adjoint group of $G$.
\end{fact}
\begin{fact}
If the theorem holds for the adjoint group of $G$, then it holds for 
$G$.\end{fact}
\begin{fact}
The adjoint group of $G$ factors as a product of simple groups. 
\end{fact}
\begin{fact}
The theorem holds for $G_1 \times G_2$ iff it holds for $G_1$ and 
$G_2$.\end{fact}
These facts show that if the theorem holds for
$G=GL_n$, $G=GSp_{2n}$, $G=SO_{2n+1}$,
and $G=SO_{2n}/ \{ \pm 1 \}$, then it will also hold for any group 
whose adjoint group
can be written as a product of the adjoint groups of these groups; 
therefore the theorem will hold for all
split connected reductive groups which are weakly classical.

The strategy of the proof of the theorem is as follows.
For each of the above groups $G$, we modify $\nu$ by a certain 
procedure to obtain
a coweight $\eta \in W\nu$.  We then show that $\eta$ is $G$-dominant
and moreover we find a Levi subgroup $L \supseteq M$ for which $\eta$
satisfies all of the hypotheses on $\nu$.
(In fact, $\eta$ satisfies all of the hypotheses on $\nu$ for $M$,
but it is easier to work with $L$.)
Thus the problem is reduced to
proving the theorem with the additional hypothesis that $\nu$ is $G$-dominant.
 We complete the proof by showing that $\nu\leq\mu$.

\section{A useful lemma}

Let $G$, $B$, $T$, and $\mathfrak{a}$ be as above.  Let $P$ be a parabolic 
subgroup of $G$ containing $B$.  Let $N$ be the unipotent radical of $P$ and 
let $M$ be the unique Levi subgroup of $P$ containing $T$; thus $P=MN$.  
As usual, let $W_M$ be the Weyl group for $(M,T)$. Let 
$\{ \alpha_i \}_{i\in I}$ be the set  of simple roots for $G$.  
Then we can write $I = I_M \sqcup I_N$,
where $I_M$, resp.~$I_N$, is the set of indices for the simple roots
 in $M$, resp.~$N$.

\begin{lem}\label{beta}
Let $\mu \in \mathfrak{a}$ be 
$G$-dominant and let $\beta\in \mathfrak{a}^{W_M}$.  
Assume further that $\mu - \beta$ is in the linear span of the coroots of 
$G$.  If  $\langle\beta, \omega_i \rangle \leq \langle \mu,\omega_i \rangle$ 
for  all fundamental weights $\omega_i$ such that $i \in I_N$, then 
$\beta\leq\mu$.
\end{lem}

\begin{proof}
Say 
\begin{eqnarray*}
\mu-\beta & = & \sum_{i\in I}r_i\alpha_i^\vee.
\end{eqnarray*}
Since $I = I_M \sqcup I_N$, we can rewrite this as 
\begin{eqnarray}
\mu-\beta & = & \sum_{i\in I_M} r_i \alpha_i^\vee  + \sum_{i\in I_N}r_i 
\alpha_i^\vee  .\label{MN equation}
\end{eqnarray}
By hypothesis, we  have $r_i \geq 0$ for all $i \in I_N$.
It remains to show that $r_i \geq 0$ for all $i \in I_M$.
To do so, it is enough to show that 
$\sum_{i\in I_M}r_i\alpha_i^\vee $ is $M$-dominant 
(cf.~\cite{Bour81} Ch.~VI \S 1, 10 ), i.e. that 
$\langle\sum_{i\in I_M}r_i\alpha_i^\vee , \alpha_j \rangle\geq 0$ 
for all $j \in I_M$.  By equation (\ref{MN equation}), this is equivalent
to showing that
\begin{eqnarray}\label{AA equation}
\langle\mu -\beta -\sum_{i\in I_N}r_i\alpha_i^\vee,\alpha_j\rangle & \geq & 0.
\end{eqnarray}
For all $j \in I_M$, we have $\langle\mu,\alpha_j\rangle\geq 0$ since $\mu$ is
 $G$-dominant, and we have  $\langle\beta, \alpha_j\rangle = 0$ since 
$\beta \in \mathfrak{a}^{W_M}$.
Finally, $\langle\sum_{i\in I_N}r_i \alpha_i^\vee, \alpha_j \rangle \leq 0$ 
since $\langle\alpha_i^\vee, \alpha_j\rangle\leq 0$ and $r_i \geq 0$ for 
all $i \in I_N$ and $j\in I_M$.  Therefore inequality~(\ref{AA equation}) 
holds for all $j\in I_M$. 
\end{proof}

\section{The proof for $G=GL_n$}

Let $G= GL_n$ and let $M$ be the Levi subgroup 
$GL_{n_1} \times  GL_{n_2} \times    \cdots \times GL_{n_r}$ of $G$ where 
$n_1 + n_2 + \cdots + n_r=n$.   Let $T$ be diagonal matrices and $B$ upper
 triangular matrices.  Then $X ={\mathbb Z}^n$.   In this case, $x\in X$ is 
$G$-dominant if $x_1\geq x_2\geq\cdots\geq x_n$.  Thus the $G$-dominant, 
$G$-minuscule elements of $X$ are of the form
$$(\underbrace{1, \ldots, 1}_t, \underbrace{0, \ldots, 0}_{n-t}) + k(1,\ldots, 1)$$
where $0 \leq t \leq n$ and $k$ is an integer.
Also $x\leq\mu$ if the following conditions hold for 
$S_i(x)=x_1 + x_2 + \cdots + x_i$:
\begin{eqnarray}
S_i(x) &\leq & S_i(\mu) \label{* equation} \quad\textrm{for } 1\leq i <n,\\
S_n(x)& =& S_n(\mu) \label{** equation}.
\end{eqnarray}

The vector $\nM$ is obtained from $\nu$ by averaging the entries of $\nu$ 
over batches where the first $n_1$ entries of $\nu$ constitute the first
 batch, the next $n_2$ the second batch, and so on.  Thus we have 
$$\nM = (\underbrace{\no,\ldots,\no}_{n_1},\underbrace{\nt,\ldots,\nt}_{n_2},
\ldots,\underbrace{\nr,\ldots,\nr}_{n_r}),$$ where $ \nik $ denotes the 
average of the entries of the $k^{\textrm{th}}$ batch of $\nu$.  We define 
$\Nk=n_1+ n_2 + \cdots + n_k$.

To prove the theorem, we will reorder the entries of $\nu$ to form a new 
coweight $\eta$ and show that, for the proper choice of Levi subgroup, $\eta$ 
satisfies all of the hypotheses on $\nu$ as well as being $G$-dominant.  This 
reduces the problem to proving the theorem with the additional hypothesis that
 $\nu$ is $G$-dominant.
We first show that the batches of $\nu$ satisfy a nice order property.

\begin{lem} \label{orderA}
Let $f_k(\nu)$ denote the first entry of the $k^{\textrm{th}}$ batch of $\nu$.
  Then $f_1(\nu)\geq f_2(\nu) \geq \ldots \geq f_r(\nu)$.
\end{lem}

\begin{proof}
 Suppose that there exists $k<r$ such that $f_{k+1}(\nu) > f_k(\nu)$.   
Both are integers, so 
\begin{eqnarray}
 f_{k+1}(\nu)-1\geq f_{k}(\nu).\label{fk equation}
\end{eqnarray} 
 Since $\nu$ is $M$-minuscule and $M$-dominant and $\nk$ is the average of the
 $k^{\textrm{th}}$ batch of $\nu$, we have 
$f_k(\nu) \geq \nk > f_k(\nu)-1$. 
 Similarly, 
\begin{equation*}
f_{k+1}(\nu) \geq \nkk > f_{k+1}(\nu)-1.
\end{equation*}
  Combining these with (\ref{fk equation}) gives 
$\nkk > f_{k+1}(\nu)-1\geq f_{k}(\nu)\geq \nk$, which contradicts the 
$G$-dominance of $\nM$. 
\end{proof}

We will now create the new coweight $\eta$ by reordering the entries of $\nu$ 
in such a way that the inequalities in Lemma~\ref{orderA} are strict for 
$\eta$.  We will then show that $\eta$ still satisfies all of the hypothesis 
on $\nu$, but for a different Levi subgroup.  To simplify doing so, we first 
prove the following lemma.

\begin{lem}\label{betaA}
 Let $\beta\in {\mathfrak a}^{W_M}$ be of the form 
$$\beta= (\underbrace{\beta_1,\ldots,\beta_1}_{n_1},\underbrace{\beta_2,
\ldots,\beta_2}_{n_2},\ldots,\underbrace{\beta_r,\ldots,\beta_r}_{n_r}).$$
To show that $\beta\leq\mu$, it is enough to show that the inequalities 
corresponding to the end of each batch are satisfied, i.e. that 
inequality~(\ref{* equation}) holds for $i=\Nk$ for all $k<r$ and that 
condition~(\ref{** equation}) holds.
\end{lem}

\begin{proof}
Follows from Lemma~\ref{beta}.
\end{proof}

Now we form the coweight $\eta$ by combining all batches of $\nu$ which have 
the same first entry into one batch and reordering each new batch in 
nonincreasing order.  Let $L$ be the Levi subgroup corresponding to these 
new batches.
Let $\eta_L$ be the vector obtained by averaging $\eta$ 
over the batches of $L$ and denote its entries by $\bar\eta_k$.  We now check 
that $\eta$ is $G$-dominant and  satisfies all of the hypotheses on $\nu$, 
but for the Levi subgroup~$L$.

\begin{lem} \label{etaA}
The coweight $\eta$ is $L$-minuscule and $\eta_L\leq\mu$.  Moreover, $\eta$ 
is  $G$-dominant, therefore $\eta$ is $L$-dominant and $\eta_L$ is 
$G$-dominant.
\end{lem}

\begin{proof}
By construction, $\eta$ is $L$-minuscule since $\nu$ is $M$-minuscule.  Also, 
the entries of $\eta$ are nonincreasing, so $\eta$ is $G$-dominant.  Therefore
$\eta$ is $L$-dominant and $\eta_L$ is $G$-dominant. 

Finally, for every $k$ there exists a $j$ so that the sum of the first $k$
 batches of $\eta$ is equal to the sum of the first $j$ batches of $\nu$.  
Therefore, since $\nM\leq\mu$, the inequalities for $\eta_L\leq\mu$ are 
satisfied at the end of each batch and by Lemma~\ref{betaA}, we have 
$\eta_L\leq\mu$.  
\end{proof}

We have shown that $\eta$ satisfies all of the hypotheses on $\nu$ for the 
Levi subgroup~$L$ and that $\eta$ is $G$-dominant.
Moreover, by its construction, $\eta \in W\nu$ so it is enough to prove 
the theorem for $(L,\eta)$ instead of $(M, \nu)$.  Thus it is enough 
to prove the theorem with the additional hypothesis that $\nu$ is 
$G$-dominant.  We can now prove that $\nu\in \textrm{Conv}(W\mu)$ by 
proving that $\nu\leq\mu$.

\begin{thm}\label{MTA}
 $\nu\in \textup{Conv}(W\mu)$
\end{thm}

\begin{proof}
We  will suppose  that $\nu \nleq \mu$ and obtain a contradiction.  If 
$\nu \nleq \mu$, then there exists an $i$ such that 
\begin{eqnarray}
\nu_1 + \nu_2 + \cdots + \nu_i &>& \mu_1 + \mu_2 + \cdots +\mu_i. \label{d equation} 
\end{eqnarray} 
Choose the smallest such $i$.  Then $\nu_i > \mu_i$ and both are integers, so 
$\nu_i -1 \geq \mu_i$.   

Suppose $\nu_i$ is in the $k^{\textrm{th}}$ batch of $\nu$.  We consider the 
$(i+1)^{\textrm{th}}$ to $\Nk^{\textrm{th}}$ entries of $\nu$ and $\mu$. 
 Since $\nu$ is $M$-dominant and $M$-minuscule, 
$\nu_{i+1}, \ldots, \nu_{\Nk} \in \{\nu_i, \nu_i -1  \}$.  
Thus 
$\nu_{i+1}+ \cdots + \nu_{\Nk} \geq  (\Nk-i)(\nu_i-1).$  
Also, since $\mu$ is $G$-dominant and $\mu_i \leq \nu_i -1$, it follows that 
$\mu_{i+1} + \cdots + \mu_{\Nk}  \leq  (\Nk-i)(\nu_i -1)$.  
Thus 
\begin{eqnarray*}
\mu_{i+1} + \cdots + \mu_{\Nk} &\leq & \nu_{i+1}+ \cdots + \nu_{\Nk}.
\end{eqnarray*}  
Combining this with inequality~(\ref{d equation}) yields 
\begin{eqnarray*}
\mu_1 + \cdots + \mu_{\Nk} < \nu_1 + \cdots + \nu_{\Nk},
\end{eqnarray*}
 which contradicts the hypothesis that $\nM \leq \mu$ since 
$\nu_1 + \cdots + \nu_{\Nk}=n_1\no + \cdots + n_r\nr$.  
 \end{proof}

\section{The proof for $G=SO_{2n+1}$}

Let $G= SO_{2n+1}$ and let $M$ be the Levi subgroup 
$GL_{n_1} \times  GL_{n_2} \times    \cdots \times GL_{n_r} \times SO_{2j+1}$ 
of $G$ where $n_1 + n_2 + \cdots + n_r+j=n$. 
Let $T$ be diagonal matrices and $B$ upper triangular matrices.  Then 
$$X =\{(a_1,a_2,\ldots,a_n,0,-a_n,\ldots,-a_2,-a_1): a_i \in {\mathbb Z}\}
\cong {\mathbb Z}^n.$$   
In this case, $x\in X$ is $G$-dominant if 
$x_1\geq x_2\geq\cdots\geq x_n \geq 0$.  Thus the $G$-dominant, $G$-minuscule 
elements of $X$ are $(1, 0, \ldots, 0)$ and $(0, \ldots, 0)$.   Also 
$x\leq\mu$ if the following condition holds for 
$S_i(x)=x_1 + x_2 + \cdots + x_i$:
\begin{eqnarray}
S_i(x) &\leq & S_i(\mu) \label{o* equation} \quad\textrm{for } 1\leq i \leq n.
\end{eqnarray}
The hypothesis that $\varphi_G(x)=\varphi_G(\mu)$ is equivalent to 
\begin{eqnarray}
S_n(\mu) - S_n(x) & \in & 2{\mathbb Z}. \label{o** equation}
\end{eqnarray}

The vector $\nM$ is obtained from $\nu$ by averaging the entries of $\nu$ over
 batches where the first $n_1$ entries of $\nu$ constitute the first batch, 
the next $n_2$ the second batch, and so on until the final batch.  The value 
for the entries of the final batch of $\nM$ is obtained by averaging over the
 middle $2j+1$ entries of $\nu$.   Thus we have 
$$\nM = (\underbrace{\no,\ldots,\no}_{n_1},\underbrace{\nt,\ldots,\nt}_{n_2},
\ldots,\underbrace{\nr,\ldots,\nr}_{n_r},\underbrace{0,\ldots,0}_j),$$ 
where $ \nik $ denotes the average of the entries of the $k^{\textrm{th}}$ 
batch of $\nu$. We define $\Nk=n_1+ n_2 + \cdots + n_k$ for $k\leq r$ 
and $\sigma(r+1) = n$.

To prove the theorem, we will reorder the entries of $\nu$ to form a new 
coweight $\eta$ and show that, for the proper choice of Levi subgroup, $\eta$ 
satisfies all of the hypotheses on $\nu$ as well as being $G$-dominant.  This 
reduces the problem to proving the theorem with the additional hypothesis 
that $\nu$ is $G$-dominant.  We first show that all of the entries of $\nu$ 
are nonnegative.

\begin{lem}\label{posB}
$\nu_i\geq 0$ for all $i$ 
\end{lem}

\begin{proof}
Since $\nu$ is $M$-minuscule and $M$-dominant, $\nu_i \geq 0 $ for $i > n-j$.
  Suppose $\nu_i <0$ for some $i\leq n-j$ and that $\nu_i$ is in the 
$k^{\textrm{th}}$ batch of $\nu$.   For $\nu_m$ in the $k^{\textrm{th}}$ 
batch, we have $\nu_m\leq \nu_i +1$  since $\nu$ is $M$-minuscule, and 
therefore $\nu_m\leq 0$ since $\nu_i$ is an negative integer.  
Thus $\bar\nu_k < 0$.  This contradicts the $G$-dominance of $\nu_M$. 
\end{proof}

Next we show that the batches of $\nu$ satisfy a nice order property.

\begin{lem}\label{orderB}
Let $f_k(\nu)$ denote the first entry of the $k^{\textrm{th}}$ batch of $\nu$.
   Then $f_1(\nu)\geq f_2(\nu)\geq\ldots\geq f_r(\nu).$
\end{lem}

\begin{proof}
The proof proceeds exactly as in Lemma~\ref{orderA}.  
\end{proof}
  
We will now create a new coweight $\eta^\prime$ by reordering the entries of 
$\nu$ in such a way that the inequalities in Lemma~\ref{orderB} are strict for
 $\eta^\prime$.  We will then modify $\eta^\prime$ slightly to form the 
coweight $\eta$ and show that $\eta$ still satisfies all of the hypotheses on 
$\nu$, but for a different Levi subgroup.  To simplify doing so, we first 
prove the following lemma.

\begin{lem}\label{betaB}
Let $\beta\in {\mathfrak a}^{W_M}$ be of the form \begin{eqnarray*}
\beta= (\underbrace{\beta_1,\ldots,\beta_1}_{n_1},\underbrace{\beta_2, 
\ldots,\beta_2}_{n_2},\ldots,\underbrace{\beta_r,\ldots,\beta_r}_{n_r}, 
\underbrace{\beta_{r+1},\ldots,\beta_{r+1}}_j),
\end{eqnarray*}
where $\beta_{r+1}=0$.  To show that $\beta\leq\mu$, it is enough to show 
that the inequalities corresponding to the end of each batch are satisfied, 
i.e. that inequality~(\ref{o* equation}) holds for $i=\Nk$ for all $k\leq r$.
\end{lem}

\begin{proof}
Follows from Lemma~\ref{beta}.
\end{proof}

Now we form the coweight $\eta$ in two steps.  First, we form the coweight 
$\eta^\prime$ by considering the first $r$ batches of $\nu$; we combine all
 batches which have the same first entry into one batch and reorder each new 
batch in nonincreasing order.  We take the final batch of $\nu$ as the final
 batch of $\eta^\prime$.  Let 
$L = GL_{m_1} \times \cdots \times GL_{m_s} \times SO_{2j+1}$ be the Levi 
subgroup corresponding to these new batches.  By construction, $\eta^\prime$ 
is $L$-dominant, and is $L$-minuscule since $\nu$ is $M$-minuscule.  Moreover 
$\eta^\prime$ is $GL_{n-j}\times SO_{2j+1}$-dominant as in Lemma~\ref{etaA}.
To form $\eta$ we modify $\eta^\prime$ by considering its final batch. 
There are two cases.  

First, if the final batch is of the form $1,0, \ldots, 0$ (so, in particular, 
$j\not= 0$) and the last entry of the $s^{\textrm{th}}$ batch of $\eta^\prime$
 is zero, then we form $\eta$ by swapping the one at the beginning of the 
final batch of $\eta^\prime$ with the left most zero entry in $\eta^\prime$.  
For example, let $G=SO_{13}$, let 
$M=GL_2 \times GL_1 \times GL_1 \times SO_5$, and let $\nu=(2,1,2,0,1,0)$.  
Then $\eta^\prime = (2,2,1,0,1,0)$, the Levi subgroup 
$L=GL_3 \times GL_1 \times SO_5$, and $\eta = (2,2,1,1,0,0)$.
Note that the one will move into a batch consisting only of zeros and ones 
since $\eta^\prime$ is $L$-minuscule and has no negative entries; hence $\eta$
 is $L$-minuscule.   

Otherwise, we set $\eta = \eta^\prime$.

We obtain $\eta_L$, resp.~$\eta^\prime_L$, from $\eta$, resp.~$\eta^\prime$,
in the same manner in which we obtained $\nM$ from $\nu$, and denote its
 entries by $\bar\eta_k$, resp.~$\bar\eta^\prime_k$.  We now check that 
$\eta$ is $G$-dominant and satisfies all of the hypotheses on $\nu$, but for
 the Levi subgroup~$L$.

\begin{lem}\label{etaB}
The coweight $\eta$ is $L$-minuscule and $\eta_L\leq\mu$.  Moreover $\eta$ 
is $G$-dominant, hence $\eta$ is $L$-dominant and $\eta_L\leq\mu$.
\end{lem}

\begin{proof}
We have already shown that $\eta$ is $L$-minuscule.  
By construction, the entries of $\eta$ are nonincreasing and, by 
Lemma~\ref{posB}, they are all nonnegative, so $\eta$ is $G$-dominant.
Therefore $\eta$ is $L$-dominant and $\eta_L$ is $G$-dominant.

It remains to show that $\eta_L \leq \mu$.  By Lemma~\ref{betaB}, it is enough
 to check that inequality~(\ref{o* equation}) holds for 
$i=\tilde\sigma(k)$ for all $k$, where $\tilde\sigma(k)=m_1 + \cdots + m_k$ 
for $k \leq s$ and $\tilde\sigma(s+1)=n$.  
We see that $\eta^\prime _L \leq \mu$ as Lemma~\ref{etaA}; therefore   
if $\eta = \eta^\prime$, we have $\eta_L \leq \mu$.
Otherwise, suppose that the left most zero entry of $\eta^\prime$ is in its 
$l^{\textrm{th}}$ batch, i.e. that we swapped the one from the final batch 
with a zero from the $l^{\textrm{th}}$ batch.  It follows that inequality 
(\ref{o* equation}) holds for $k < l$ as Lemma~\ref{etaA}.  Moreover, since 
$\eta$ is $G$-dominant and all of its entries are nonnegative, all of the 
entries to the right of the one must be zero.  Hence, since $\mu_i\geq 0$ 
for all $i$ as $\mu$ is $G$-dominant, it is enough to check inequality 
(\ref{o* equation}) for $i=\tilde\sigma(l)$.

Since $\eta^\prime_L\leq\mu$, we have 
$S_{\tilde\sigma(l)}(\eta^\prime_L)\leq S_{\tilde\sigma(l)}(\mu)$.  
If $S_{\tilde\sigma(l)}(\eta^\prime_L)< S_{\tilde\sigma(l)}(\mu)$, 
then $S_{\tilde\sigma(l)}(\eta_L) \leq S_{\tilde\sigma(l)}(\mu)$ since both 
sides are integral and  
$S_{\tilde\sigma(l)}(\eta_L) = S_{\tilde\sigma(l)}(\eta^\prime_L) + 1$.  
Otherwise, 
\begin{eqnarray}
S_{\tilde\sigma(l)}(\eta^\prime_L) &=& S_{\tilde\sigma(l)}(\mu). \label{ot equation}
\end{eqnarray}
Since $\eta^\prime_L\leq\mu$, we have
\begin{eqnarray}
S_{\tilde\sigma(l)-1}(\eta^\prime_L) &\leq& S_{\tilde\sigma(l)-1}(\mu).\label{ox equation}  
\end{eqnarray}
Combining this with equality (\ref{ot equation}) yields 
\begin{eqnarray}
\bar\eta^\prime_l &\geq& \mu_{\tilde\sigma(l)}.\label{ou equation}
\end{eqnarray}
Since the $l^{\textrm{th}}$ batch of $\eta^\prime$ consists only of zeros 
and ones and contains at least one zero, we have $\bar\eta^\prime _l < 1$.  
We will obtain a contradiction by showing that $\mu_{\tilde\sigma(l)}\geq 1$.
Now, since $\eta^\prime_L\leq\mu$, we have
\begin{eqnarray}
S_{\tilde\sigma(s)}(\eta^\prime_L) + 1 &\leq& S_{\tilde\sigma(s)+1}(\mu). \label{ow equation}
\end{eqnarray}  
Therefore $1 \leq \mu_{\tilde\sigma(l)+1} + \cdots + \mu_{\tilde\sigma(s)+1}$ 
by inequalities (\ref{ot equation}) and (\ref{ow equation}) and since
$\bar\eta^\prime_k=0$ for all $k>l$.
Hence, since $\mu$ is $G$-dominant, it follows that 
$\mu_{\tilde\sigma(l)}\geq\mu_{\tilde\sigma(l)+1}\geq 1$, giving the 
desired contradiction.  Thus 
$S_{\tilde\sigma(l)}(\eta^\prime_L)< S_{\tilde\sigma(l)}(\mu)$ and 
$\eta_L\leq\mu$.
\end{proof}

We have shown that $\eta$ satisfies all of the hypotheses on $\nu$ for the 
Levi subgroup~$L$ and that $\eta$ is $G$-dominant.
Moreover, by its construction, $\eta \in W\nu$ so it is enough to prove the 
theorem for $(L,\eta)$ instead of $(M, \nu)$.  Thus it is enough to prove the 
theorem with the additional hypothesis that $\nu$ is $G$-dominant.
We can now prove that $\nu\in \textrm{Conv}(W\mu)$ by proving that 
$\nu\leq\mu$.

\begin{thm}\label{MTB}
$\nu\in \textup{Conv}(W\mu)$
\end{thm}

\begin{proof}
We will suppose $\nu \nleq \mu$ and obtain a contradiction.  
If $\nu \nleq \mu$, then there exists an $i$ such that 
\begin{eqnarray}
\nu_1 + \nu_2 + \cdots + \nu_i &>& \mu_1 + \mu_2 + \cdots +\mu_i. 
\end{eqnarray} 
Choose the smallest such $i$ and suppose that $\nu_i$ is in the 
$k^{\textrm{th}}$ batch of $\nu$.  As in Theorem~\ref{MTA}, we obtain 

\begin{eqnarray}
\mu_1+\cdots + \mu_{\Nk} &<& \nu_1 + \cdots +\nu_{\Nk},\label{or equation}  
\end{eqnarray} 
which produces a contradiction as follows:
 If $\Nk=n$ and the first entry of the final batch is one, then  
$\nu_1 + \cdots +\nu_{\Nk}=n_1\bar\nu_1+\cdots + n_k\nk + 1$.  
Thus, since both values are integers,  
$\mu_1+\cdots + \mu_{\Nk} \leq  n_1\bar\nu_1+\cdots + n_k\nk$.  
Since $\nM\leq\mu$, we must have equality.  This contradicts condition 
(\ref{o** equation}). Otherwise 
$\nu_1 + \cdots +\nu_{\Nk}=n_1\bar\nu_1+\cdots + n_k\nk$, so 
inequality~(\ref{or equation}) contradicts the hypothesis that $\nM\leq\mu$.
\end{proof}

\section{The proof for $G=SO_{2n}/\{\pm 1 \}$ (Integers)}\label{Dsect}

Let $G= SO_{2n}/\{ \pm 1 \}$ and let $M$ be the Levi subgroup 
$GL_{n_1} \times  GL_{n_2} \times    \cdots \times GL_{n_r} \times SO_{2j}$ 
of $G$ where $n_1 + n_2 + \cdots + n_r+j=n$.  We may assume that $j\not= 1$ 
since the case $j=1$ is the same as the one in which $j=0$ and $n_r=1$.
There are other Levi subgroups in $G$, but it is not necessary
 to consider them since there exists an outer automorphism of $SO_{2n}$ 
which takes each of these to a Levi subgroup which we have considered.
Let $T$ be diagonal matrices and $B$ upper triangular matrices.  Then 
$$X =\{(a_1,a_2,\ldots,a_n,-a_n,\ldots,-a_2,-a_1): \textrm{either }a_i 
\in {\mathbb Z} \ \forall\  i \textrm{ or }   a_i \in 
\frac{1}{2}({\mathbb Z} \smallsetminus 2{\mathbb Z} ) \ \forall \  i\}. $$ 
We first consider the case in which $a_i \in {\mathbb Z}$ for all $i$; 
we will consider the other case in Section~\ref{DHsect}.
In this case, $x\in X$ is $G$-dominant if 
\begin{eqnarray*}
x_1\geq x_2 \geq \cdots\geq x_{n-1} &\geq& x_n, \\
\textrm{and } x_{n-1}+x_n &\geq& 0. 
\end{eqnarray*}
It follows that if $x$ is $G$-dominant, then $x_i \geq 0$ for all $i\leq n-1$.
The $G$-dominant, $G$-minuscule elements of $X$ are $(1,0,\ldots , 0)$ and 
$(0, \ldots , 0)$.
Also $x\leq\mu$ if the following conditions hold for 
$S_i(x)=x_1 + x_2 + \cdots + x_i$:
\begin{eqnarray}
S_i(x) &\leq & S_i(\mu) \label{e* equation} \quad\textrm{for } 1\leq i \leq n-2,\\
S_{n-1}(x)-x_n &\leq & S_{n-1}(\mu)-\mu_n, \label{e- equation} \\
S_n(x) &\leq & S_n(\mu),\label{en equation}.
\end{eqnarray}
The hypothesis that $\varphi_G(x) = \varphi_G(\mu)$ is equivalent to 
\begin{eqnarray}
S_n(\mu) - S_n(x) & \in & 2{\mathbb Z} \label{e** equation}
\end{eqnarray}

The vector $\nM$ is obtained from $\nu$ by averaging the entries of $\nu$ over
 batches where the first $n_1$ entries of $\nu$ constitute the first batch, 
the next $n_2$ the second batch, and so on until the final batch.  The value 
for the entries of the final batch of $\nM$ is obtained by averaging over the 
middle $2j$ entries of $\nu$.  Thus we have 
$$\nM = (\underbrace{\no,\ldots,\no}_{n_1},\underbrace{\nt,\ldots,\nt}_{n_2},
\ldots,\underbrace{\nr,\ldots,\nr}_{n_r},\underbrace{0,\ldots,0}_j),$$ 
where $ \nik $ denotes the average of the entries of the $k^{\textrm{th}}$ 
batch of $\nu$.  We define $\Nk=n_1+ n_2 + \cdots + n_k$ for $k\leq r$ and 
$\sigma(r+1) = n$.

To prove the theorem, we will reorder the entries of $\nu$ to form a new 
coweight $\eta$ and show that, for the proper choice of Levi subgroup, $\eta$ 
satisfies all of the hypotheses on $\nu$ as well as being $G$-dominant.  This 
reduces the problem to proving the theorem with the additional hypothesis that
 $\nu$ is $G$-dominant.  We first show that at most one of the entries of 
$\nu$ is negative.

\begin{lem}\label{posD}
$\nu_i\geq 0$ for all $i$, unless $j=0$ and $n_r=1$, in which case, 
$\nu_i\geq 0$ for all $i\leq n-1$
\end{lem}

\begin{proof}
Suppose that $\nu_i < 0$ and that $\nu_i$ is in the $k^{\textrm{th}}$ batch of
 $\nu$.  As in Lemma~\ref{posB}, we see that $k \not= r+1$ and that 
$\bar\nu_k <0$.  Since $\nM$ is $G$-dominant, $\bar\nu_k < 0$ can occur only 
if $j=0$, $n_r=1$, and $k=r$, in which case $i=n$.
\end{proof}

Next we show that the batches of $\nu$ satisfy a nice order property.

\begin{lem}\label{orderD}
Let $f_k(\nu)$ denote the first entry of the 
$k^{\textrm{th}}$ batch of $\nu$.   Then 
$f_1(\nu)\geq f_2(\nu)\geq\ldots\geq f_r(\nu).$
\end{lem}

\begin{proof}
The proof proceeds exactly as Lemma~\ref{orderA}. 
\end{proof}
  
We will now create a coweight $\eta^\prime$ by reordering the entries of 
$\nu$ in such a way that the inequalities in Lemma~\ref{orderD} are strict for
 $\eta^\prime$.  We will then modify $\eta^\prime$ slightly to form the 
coweight $\eta$ and show that $\eta$ still satisfies all of the hypothesis 
on $\nu$, but for a different Levi subgroup.  To simplify doing so, we first 
prove the following lemma.

\begin{lem}\label{betaD}
Let $\beta\in {\mathfrak a}^{W_M}$ be of the form
\begin{eqnarray*}
\beta= (\underbrace{\beta_1,\ldots,\beta_1}_{n_1},
\underbrace{\beta_2,\ldots,\beta_2}_{n_2},\ldots,
\underbrace{\beta_r,\ldots,\beta_r}_{n_r}, \
\underbrace{\beta_{r+1},\ldots,\beta_{r+1}}_j),
\end{eqnarray*}
where $\beta_{r+1}=0$.  To show that $\beta\leq\mu$, it is enough to show that
 the inequalities corresponding to the end of each batch are satisfied, i.e. 
that inequality~(\ref{e* equation}) holds for $i=\Nk$ for all  $k$
 such that  $\sigma(k)\leq n-2$,  that inequality~(\ref{en equation}) holds, 
and, if $j=0$ and $n_r=1$, that inequality~(\ref{e- equation}) holds.
\end{lem}

\begin{proof}
Follows from Lemma~\ref{beta}. 
\end{proof}

Now we form the coweights $\eta^\prime$ and $\eta$ as for $G=SO_{2n+1}$.  
Again, we let $L = GL_{m_1} \times \cdots \times GL_{m_s} \times SO_{2j}$ be 
the Levi subgroup corresponding to the new batches,
We obtain $\eta_L$, resp.~$\eta^\prime_L$, from $\eta$, resp.~$\eta^\prime$, 
in the same manner in which we obtained $\nM$ from $\nu$, and denote its 
entries by $\bar\eta_k$, resp.~$\bar\eta^\prime_k$.

We now check that $\eta$ is $G$-dominant and satisfies all of the hypotheses 
on $\nu$, but for the Levi subgroup~$L$.

\begin{lem}\label{etaD}
The coweight $\eta$ is $L$-minuscule and $\eta_L\leq\mu$.  Moreover $\eta$ is 
$G$-dominant, hence $\eta$ is $L$-dominant and $\eta_L\leq\mu$.
\end{lem}

\begin{proof}
As in the $SO_{2n+1}$ case, we see that $\eta$ and $\eta^\prime$ are 
$L$-minuscule.  We claim that $\eta$ is $G$-dominant.  As before, the 
inequalities  $\eta_1 \geq \cdots \geq\eta_n$ follow from the way $\eta$ 
was constructed, so we need only check that $\eta_{n-1} + \eta_n \geq 0$.  
This is automatic unless some entry of $\nu$ is strictly negative, so by 
Lemma~\ref{posD}  we may now assume that $j=0$ and $n_r=1$.  
 Thus we have $\eta_n = \bar\nu_r$.
We know that $\bar\nu_{r-1} + \bar\nu_r \geq 0$ since $\nM$ is $G$-dominant.  
Hence  $\bar\eta_{s-1} + \eta_n \geq 0$ since
$\bar\eta_{s-1} \geq \bar\nu_{r-1}$.
Finally, since $\eta$ is $L$-minuscule, $\eta_{n-1}$ is the greatest integer 
less than or equal to $\bar\eta_{s-1}$, so $\eta_{n-1}+ \eta_n \geq 0$.
Thus $\eta$ is $G$-dominant.  Therefore $\eta$ is $L$-dominant and $\eta_L$ 
is $G$-dominant.

It only remains to show that $\eta_L\leq\mu$.  
To do so we apply Lemma~\ref{betaD}.  
Define $\tilde\sigma(k)=m_1 + \cdots + m_k$ for $k\leq s$ and
$\tilde\sigma(s+1)=n$.
The method used in Lemma~\ref{etaB} establishes
inequality~(\ref{e* equation})  for $i=\tilde\sigma(k)$ for $k$ such that  
$\tilde\sigma(k) \leq n-2$, as well as inequality~(\ref{en equation}).  It 
remains  to verify inequality~(\ref{e- equation}) under the assumption that 
$j=0$ and $m_s=1$. In this case, we have $\eta=\eta^\prime$, and we must have 
$n_r=1$.  Thus  $\bar\eta_s=\bar\nu_r$ and  
$S_{\tilde\sigma(s-1)}(\eta_L)= S_{\sigma(r-1)}(\nu_M)$.  Therefore, since 
$\nM\leq\mu$, inequality~(\ref{e- equation}) holds.  
\end{proof}

We have shown that $\eta$ satisfies all of the hypotheses on $\nu$ for the 
Levi subgroup~$L$ and that $\eta$ is $G$-dominant.
Moreover, by its construction, $\eta \in W\nu$ so it is enough to prove the 
theorem for $(L,\eta)$ instead of $(M, \nu)$.  Thus it is enough to prove the 
theorem with the additional hypothesis that $\nu$ is $G$-dominant.
We can now prove that $\nu\in \textrm{Conv}(W\mu)$ by proving that
 $\nu\leq\mu$.

\begin{thm}\label{MTD}
$\nu\in \textup{Conv}(W\mu)$
\end{thm}

\begin{proof}
We will suppose $\nu \nleq \mu$ and obtain a contradiction. 
If $\nu \nleq \mu$, then either there exists an $i \not= n-1$ such that 
\begin{eqnarray}
\nu_1 + \nu_2 + \cdots + \nu_i &>& \mu_1 + \mu_2 + \cdots +\mu_i \label{em equation}
\end{eqnarray} 
or inequality~(\ref{e- equation}) fails.

If inequality~(\ref{em equation}) holds for $i=n$, then we will argue as in 
Theorem~\ref{MTB}.
If  the first entry of the final batch is zero, then 
$\nu_1 + \cdots +\nu_{n}=n_1\bar\nu_1+\cdots + n_k\nk$, so 
inequality~(\ref{em equation}) contradicts the hypothesis that $\nM\leq\mu$.
 If the first entry of the final batch is one, then  
$\nu_1 + \cdots +\nu_{n}=n_1\bar\nu_1+\cdots + n_k\nk + 1$.  Combining this 
with  inequality~(\ref{em equation})  yields 
$\mu_1+\cdots + \mu_{n} \leq  n_1\bar\nu_1+\cdots + n_k\nk$ since both values 
are integers.  Since $\nM\leq\mu$, we must have equality. 
 This contradicts condition~(\ref{e** equation}).

Now suppose that there exists an $i\leq n-2$ for which 
inequality~(\ref{em equation}) holds.
Choose the smallest such $i$ and suppose that $\nu_i$ is in the 
$k^{\textrm{th}}$ batch of $\nu$.  As in Lemma~\ref{MTA}, we obtain 
\begin{eqnarray}
\mu_1+\cdots + \mu_{\Nk} &<& \nu_1 + \cdots +\nu_{\Nk}.\label{er equation}  
\end{eqnarray} 
If $\sigma(k) \not= n-1$, then we obtain a contradiction as for $SO_{2n+1}$.

If $\sigma(k)= n-1$, then we have $j=0$ and $n_r=1$. 
 We will first show that 
inequality~(\ref{e- equation}) holds in this case.
Since
$\nM\leq\mu$, we have
\begin{eqnarray}
S_{\sigma(r-1)}(\nM)-\nr \leq S_{n-1}(\mu)- \mu_n. \label{eA equation}
\end{eqnarray}
By the definition of $\nM$, the partial sums are equal at the end of 
each batch and we are assuming that $j=0$ and $n_r=1$, so
$\nr =\nu_n  $ and $S_{\sigma(r-1)}(\nM)=S_{n-1}(\nu)$.  Substituting these 
values into inequality~(\ref{eA equation}) gives inequality 
(\ref{e- equation}) as desired.
Adding inequality~(\ref{e- equation}) to inequality~(\ref{en equation})
 yields $S_{n-1}(\nu) \leq S_{n-1}(\mu)$, contradicting 
inequality~(\ref{er equation}). 

Now suppose that inequality~(\ref{e- equation}) fails.  Thus we have
\begin{eqnarray}
S_{n-1}(\nu)-\nu_n &>& S_{n-1}(\mu)-\mu_n.  \label{ei equation}
\end{eqnarray}
By what we have already shown, we have  
\begin{eqnarray}
S_{n-2}(\nu) &\leq & S_{n-2}(\mu). \label{eC equation}
\end{eqnarray}
We claim that
\begin{eqnarray}
\nu_{n-1}-\nu_n &>& \mu_{n-1}-\mu_n \geq 0.  \label{ej equation}
\end{eqnarray}
The first inequality follows from inequalities (\ref{ei equation}) and
(\ref{eC equation}), and the second holds since $\mu$ is $G$-dominant.
We have shown that inequality~(\ref{e- equation}) holds if $j=0$ and $n_r=1$ 
so we may assume that this is not the case.  Thus $\nu_{n-1}$ and $\nu_n$ 
are in the same batch of $\nu$, so, since $\nu$ is $M$-minuscule,
$\nu_{n-1}-\nu_n\in \{ 0, 1 \}$.
Combining this with  inequality~(\ref{ej equation}) 
gives  
\begin{eqnarray}
\nu_{n-1}- \nu_n =1  & \textrm{and} & \mu_{n-1} - \mu_n =0. \label{eD equation}\end{eqnarray}
Substituting these values into 
inequality  (\ref{ei equation}) and  combining the result with 
inequality~(\ref{eC equation})  gives
\begin{eqnarray*}
S_{n-2}(\nu)+1>S_{n-2}(\mu)\geq S_{n-2}(\nu).  
\end{eqnarray*}
Both are integers, so $S_{n-2}(\mu) =S_{n-2}(\nu)$.  This contradicts 
condition~(\ref{e** equation}) since it follows from the equalities in 
(\ref{eD equation}) that $S_{n}(\mu)$ has the same parity as $S_{n-2}(\mu)$
 and $S_{n}(\nu)$ has the opposite parity from $S_{n-2}(\nu)$.  
Thus $\nu\leq\mu$. 
\end{proof}

\section{The proof for $G=SO_{2n}/\{\pm 1\}$ (Half Integers)}\label{DHsect}

Recall from Section~\ref{Dsect} that for $G=SO_{2n}/ \{ \pm 1\}$, we have 
$$X =\{(a_1,a_2,\ldots,a_n,-a_n,\ldots,-a_2,-a_1): \textrm{either }
a_i\in {\mathbb Z}\ \forall\  i \textrm{ or }a_i\in \frac{1}{2}({\mathbb Z} 
\smallsetminus 2{\mathbb Z}) \ \forall\  i\}. $$
We now consider the case in which $ a_i \in \frac{1}{2}({\mathbb Z}
 \smallsetminus 2{\mathbb Z} )$ for all $i$.
As in Section~\ref{Dsect}, we let $M$ be the Levi subgroup $GL_{n_1} 
\times  GL_{n_2} \times    \cdots \times GL_{n_r} \times SO_{2j}$ of $G$ 
where $n_1 + n_2 + \cdots + n_r+j=n$ and assume that $j\not= 1$.
Rather than consider $X$, we will double the entries of $X$ and adjust the 
notion of minuscule and the definition of the partial order accordingly.   

First we note that all of the elements of $X$ that we are considering in 
this case are now $2n$-tuples of odd integers.  
As before, $x\in X$ is $G$-dominant if 
\begin{eqnarray*}
x_1\geq x_2 \geq \cdots\geq x_{n-1} &\geq & x_n, \\
\textrm{and } x_{n-1}+x_n &\geq& 0 
\end{eqnarray*}
It follows that if $x$ is $G$-dominant, then $x_i \geq 1$ for all $i\leq n-1$.
Since we have multiplied $X$ by two, we will now be concerned with elements 
whose pairing with any root yields two, zero, or negative two.  We will refer 
to these elements as $2$-minuscule.
The $G$-dominant, $G$-$2$-minuscule elements of $X$ are $(1,1,\ldots , 1)$ 
and $(1, \ldots 1, -1)$.
Also $x\leq\mu$ if the following conditions hold for
 $S_i(x)=x_1 + x_2 + \cdots + x_i$:
\begin{eqnarray}
S_i(x) &\leq & S_i(\mu) \label{eh* equation} \quad\textrm{for } 1\leq i \leq n-2,\\
S_{n-1}(x)-x_n &\leq & S_{n-1}(\mu)-\mu_n, \label{eh- equation} \\
S_n(x) &\leq & S_n(\mu),\label{ehn equation}. 
\end{eqnarray}
The hypothesis that $\varphi_G(x)=\varphi_G(\mu)$ is equivalent to
\begin{eqnarray}
S_n(\mu) - S_n(x) & \in & 4{\mathbb Z}. \label{eh** equation}
\end{eqnarray}
As in the previous sections, the first $n_1$ entries of $\nu$ constitute the 
first batch, the next $n_2$ the second batch, and so on.
The $M$-dominant, $M$-$2$-minuscule elements of $X$ will now be those such 
that the entries of each batch are nonincreasing and differ by either zero 
or two and 
such that  the $r+1^{\textrm{st}}$ batch is of the form $(1,1,\ldots , 1)$ 
or $(1, \ldots 1, -1)$.

The vector $\nM$ is obtained from $\nu$ by averaging the entries of $\nu$ over
 batches with the exception that the value for the entries of the final batch 
of $\nM$ is obtained by averaging over the middle $2j$ entries of $\nu$.  Thus
 we have 
$$\nM = (\underbrace{\no,\ldots,\no}_{n_1},\underbrace{\nt,\ldots,\nt}_{n_2},
\ldots,\underbrace{\nr,\ldots,\nr}_{n_r},\underbrace{0,\ldots,0}_j),$$ 
where $ \nik $ denotes the average of the entries of the $k^{\textrm{th}}$ 
batch of $\nu$. 
We define $\Nk=n_1+ n_2 + \cdots + n_k$ for $k\leq r$ and $\sigma(r+1) = n$.

To prove the theorem, we will reorder the entries of $\nu$ to form a new 
coweight $\eta$ and show that, for the proper choice of Levi subgroup,
 $\eta$ satisfies all of the hypotheses on $\nu$ as well as being 
$G$-dominant.  This reduces the problem to proving the theorem with the 
additional hypothesis that $\nu$ is $G$-dominant.  We first show that at 
most one of the entries of $\nu$ is less than negative one.

\begin{lem}\label{neg1D}
$\nu_i\geq -1$ for all $i$, unless $j=0$ and $n_r=1$, in which case, 
$\nu_i\geq -1$ for all $i\leq n-1$
\end{lem}

\begin{proof}
Since $\nu$ is $M$-$2$-minuscule and $M$-dominant, we have $\nu_i \geq -1 $
 for $i > n-j$. 
 Suppose that $\nu_i <-1$ for some $i\leq n-j$ and that $\nu_i$ is in the
 $k^{\textrm{th}}$ batch of $\nu$.
For $\nu_m$ in the $k^{\textrm{th}}$ batch, we have $\nu_m\leq \nu_i +2$  
since $\nu$ is $M$-$2$-minuscule, and therefore
$\nu_m\leq -1$ since $\nu_i <-1$ is an odd integer.  
Thus $\bar\nu_k <-1$. Since $\nM$ is $G$-dominant,  $\bar\nu_k < 0$ 
can only occur if $j=0$,  $n_r=1$, and $k=r$, in which case $i=n$.
\end{proof}

Next we show that the batches of $\nu$ satisfy a nice order property.

\begin{lem}\label{orderD*}
Let $f_k(\nu)$ denote the first entry of the $k^{\textrm{th}}$ batch of $\nu$.
  Then $f_1(\nu)\geq f_2(\nu) \geq \ldots \geq f_r(\nu)$.
\end{lem}

\begin{proof}
Suppose that there exists $k<r$ such that $f_{k+1}(\nu) > f_k(\nu)$.   
Both are odd integers, so
\begin{eqnarray}
f_{k+1}(\nu)-2\geq f_{k}(\nu)\label{ehfk equation}
\end{eqnarray}
Since $\nu$ is $M$-$2$-minuscule and $M$-dominant and $\nk$ is the average 
of the $k^{\textrm{th}}$ batch of $\nu$, we have 
$f_k(\nu) \geq \nk > f_k(\nu)-2$.  Similarly, 
$f_{k+1}(\nu) \geq \nkk > f_{k+1}(\nu)-2$.
  Combining these with inequality~(\ref{ehfk equation}) gives 
$\nkk > f_{k+1}(\nu)-2\geq f_{k}(\nu)\geq \nk$, which contradicts the
 $G$-dominance of $\nM$.
\end{proof}

We will now create a coweight $\eta^\prime$ by reordering the entries of 
$\nu$ in such a way that the inequalities in Lemma~\ref{orderD*} are strict 
for $\eta^\prime$.  We will then modify $\eta^\prime$ slightly to form the
 coweight $\eta$ and show that $\eta$ still satisfies all of the hypothesis 
on $\nu$, but for a different Levi subgroup.  To simplify doing so, we first 
prove the following lemma.

\begin{lem}\label{betaD*}
Let $\beta\in {\mathfrak a}^{W_M}$ be of the form
\begin{eqnarray*}
\beta= (\underbrace{\beta_1,\ldots,\beta_1}_{n_1},\underbrace{\beta_2,\ldots,
\beta_2}_{n_2},\ldots,\underbrace{\beta_r,\ldots,\beta_r}_{n_r}, 
\underbrace{\beta_{r+1},\ldots,\beta_{r+1}}_j),
\end{eqnarray*}
where $\beta_{r+1}=0$.  To show that $\beta\leq\mu$, it is enough to show 
that the inequalities corresponding to the end of each batch are satisfied, 
i.e. that inequality~(\ref{eh* equation}) holds for $i=\Nk$ for all 
$\sigma(k)\leq n-2$, that inequality~(\ref{ehn equation}) holds, and, if 
$j=0$ and $n_r=1$, that inequality~(\ref{eh- equation}) holds.
\end{lem}

\begin{proof}
Follows from Lemma~\ref{beta}. 
\end{proof}

Now we form the coweight $\eta$ in three steps.  First, we form the coweight
$\eta^\prime$ as in the $SO_{2n+1}$ case. 
 Again, we let $L = GL_{m_1} \times \cdots \times GL_{m_s} \times SO_{2j}$ be 
the Levi subgroup corresponding to the new batches.
By construction, $\eta^\prime$ is $L$-$2$-minuscule
since $\nu$ is $M$-$2$-minuscule.
It follows from Lemma~\ref{neg1D} and the construction of $\eta^\prime$ 
that $\eta^\prime_i \geq -1$ for all $i \leq n-1$.

Moreover we claim that any negative ones in $\eta^\prime$ will be in its 
 $s^{\textrm{th}}$  or $s+1^{\textrm{st}}$ batch.   
Suppose that $\eta^\prime$ contains a negative one prior to the 
$s^{\textrm{th}}$ batch.   Since the entries of $\eta^\prime$ are 
nonincreasing, this implies that all  of the entries of the 
$s^{\textrm{th}}$ batch will be at most negative one.   
The  $s^{\textrm{th}}$ batch of $\eta^\prime$ was formed by combining batches
 of $\nu$,  so this yields that  all of the entries of the 
$r^{\textrm{th}}$ batch of $\nu$ will also be at most negative one. 
Thus $\bar\nu_r < 0$ which contradicts the 
$G$-dominance of $\nM$ unless $j=0$ and $n_r=1$.  In this  case,  we have 
$\bar\nu_r = \nu_n$ and since $\nM$ is $G$-dominant, 
$\bar\nu_{r-1}+\bar\nu_r \geq 0$.  
Therefore $\bar\nu_{r-1}\geq 1$, so since $\nM$ is $G$-dominant, 
$\bar\nu_{k}\geq 1$ for all $k\leq r-1$.
Thus, since $\nu$ is $M$-$2$-minuscule, $\nu_i \geq 1$ for all $i\leq n-1$ and
it follows from the construction of $\eta^\prime$ that  $\eta^\prime_i \geq 1$
 for all $i\leq n-1$.
Therefore $\eta^\prime$ cannot contain a negative one prior to the 
$s^{\textrm{th}}$ batch.
 
Next we form a coweight $\eta^{\prime\prime}$ by replacing every negative one 
in 
the $s^{\textrm{th}}$ batch of $\eta^\prime$ with a positive
 one.  Finally, if we made an odd number of sign changes in the previous step,
 then we  change 
the sign of the final entry of $\eta^{\prime\prime}$ to form $\eta$; 
otherwise, $\eta = \eta^{\prime\prime}$.

We obtain $\eta_L$, resp.~$\eta^\prime_L$, from $\eta$, resp.~$\eta^\prime$, 
in the same manner as we obtained $\nM$ from $\nu$,
and denote its entries by $\bar\eta_k$, resp.~$\bar\eta^\prime_k$.
We now check that $\eta$  is $G$-dominant and satisfies all of the hypotheses 
on $\nu$, but for the Levi subgroup~$L$.

\begin{lem}
The coweight $\eta$ is $L$-$2$-minuscule and $\eta_L\leq\mu$.
 Moreover $\eta$ is $G$-dominant, hence $\eta$ is $L$-dominant and
$\eta_L$ is $G$-dominant. 
\end{lem}

\begin{proof}
We begin by showing that $\eta$ is $L$-$2$-minuscule and $G$-dominant. 
It follows that $\eta$ is $L$-dominant and
$\eta_L$ is $G$-dominant.

First, if $j=0$ and $n_r=1$, then $\eta = \eta^\prime$, so we have shown that 
$\eta$ is $L$-$2$-minuscule.  We claim that $\eta$ is $G$-dominant.
The inequalities 
$\eta_1 \geq \cdots \geq\eta_n$ follow from the way $\eta$ was constructed, 
so we need only check  that $\eta_{n-1} + \eta_n \geq 0$.  
As in Section~\ref{Dsect}, we obtain $\bar\eta_{s-1} + \eta_n \geq 0$.  
Since $\eta$ is $L$-$2$-minuscule,
$\eta_{n-1}$ is the greatest odd integer less than or equal to 
$\bar\eta_{s-1}$;  therefore we have
$\eta_{n-1}+ \eta_n \geq 0$ since $\eta_n$ is an odd integer.

Otherwise, it is clear that the $s+1^{\textrm{st}}$ batch of $\eta$ is still 
$L$-$2$-minuscule.  If the  $s^{\textrm{th}}$ batch of $\eta^\prime$ 
contains any negative ones, then,  since we have shown that
 $\eta^\prime$ is $L$-$2$-minuscule and by Lemma~\ref{neg1D}, all other 
entries of the batch must be 
either positive or negative one.
Therefore, since we only  changed the signs of  positive and negative ones, 
$\eta$ is $L$-$2$-minuscule.  Moreover, by construction, $\eta_i \geq 1$ for 
all $i<n$, the entries of $\eta$ are nonincreasing, and, by 
Lemma~\ref{neg1D}, $\eta_n \geq -1$; hence it is clear that $\eta$ is 
$G$-dominant.

It only remains to show that $\eta_L\leq\mu$.  To do so, we apply
 Lemma~\ref{betaD*}.  Define $\tilde\sigma(k)=m_1 + \cdots + m_k$
for $k \leq s$ and $\tilde\sigma(s+1)=n$.

First, we show that $\eta^\prime _L \leq \mu$.  The method used in 
Lemma~\ref{etaA} establishes
inequality~(\ref{eh* equation})  for $i=\tilde\sigma(k)$ for $k$ such that  
$\tilde\sigma(k) \leq n-2$, as well as inequality~(\ref{ehn equation}).  It 
remains to verify inequality~(\ref{eh- equation}) under the assumption 
that $j=0$ and $m_s=1$; the proof proceeds exactly as in Lemma~\ref{etaD}.
Thus if $\eta = \eta^\prime$, then we have $\eta_L \leq \mu$.

Now we consider the case in which $\eta \not= \eta^\prime$.  (In particular, 
it is not the case that $j=0$ and $m_s=1$.)
Since $\eta$ and $\eta^\prime$ do not differ prior to the $s^{\textrm{th}}$ 
batch and $\eta^\prime _L \leq \mu$, inequality~(\ref{eh* equation}) is
satisfied for $i=\tilde\sigma(k)$ for all $k < s$.  Moreover, since either the
 $s^{\textrm{th}}$ batch is the final batch (if $j=0$) or $\bar\eta_{s+1}=0$ 
(if $j\not= 0$), and since $\mu$ is $G$-dominant, it is enough to check that 
inequality~(\ref{eh* equation}) holds for 
$i=\tilde\sigma(s)$.

We have shown that 
$S_{\tilde\sigma(s-1)}(\eta_L) \leq S_{\tilde\sigma(s-1)}(\mu)$.  Since 
the last entry of the $s^{\textrm{th}}$ batch of $\eta$ is either positive or 
negative one and $\eta$ is $L$-$2$-minuscule, we also have that 
$\bar\eta_s \leq 1$.  Moreover, since $\mu$ is $G$-dominant, we have 
$\mu_i \geq 1$ for all $i < n$.  It follows that 
$S_{\tilde\sigma(s)}(\eta_L) \leq S_{\tilde\sigma(s)}(\mu)$ unless 
$\tilde\sigma(s)=n$.
If $\tilde\sigma(s)=n$ (so, in particular, $j=0$), then 
$$\mu_{\tilde\sigma(s-1)+1}+\cdots +\mu_n \geq m_s -2$$ 
since $\mu$ is $G$-dominant.  Thus $S_{n}(\eta_L) \leq S_{n}(\mu)$ unless  
$S_{\tilde\sigma(s-1)}(\eta_L) = S_{\tilde\sigma(s-1)}(\mu)$, $\bar\eta_s =1$,
 and $\mu_{\tilde\sigma(s-1)+1}+\cdots +\mu_n = m_s -2$.  In this case, 
we have $S_{n}(\eta_L)=S_{n}(\mu)+2$.
This is a contradiction to condition~(\ref{eh**  equation}) since 
$S_{n}(\eta_L)=S_{n}(\eta)$ and we see that $S_n(\eta)$ and $S_n(\nu)$ are 
congruent modulo four since the difference in the two sums comes from changing
 an even number of negative ones to positive ones.
\end{proof}

We have shown that $\eta$ satisfies all of the hypotheses on $\nu$ for the 
Levi subgroup~$L$ and that $\eta$ is $G$-dominant.
Moreover, by its construction, $\eta \in W\nu$ so it is enough to prove the 
theorem for $(L,\eta)$ instead of $(M, \nu)$.  Thus it is enough to prove the 
theorem with the additional hypothesis that $\nu$ is $G$-dominant.
We can now prove that $\nu\in \textrm{Conv}(W\mu)$ by proving that 
$\nu\leq\mu$.

\begin{thm}
$\nu\in \textup{Conv}(W\mu)$
\end{thm}

\begin{proof}
We will suppose $\nu \nleq \mu$ and obtain a contradiction.  
If $\nu \nleq \mu$, then either there exists an $i \not= n-1$ such that 
\begin{eqnarray}
\nu_1 + \nu_2 + \cdots + \nu_i &>& \mu_1 + \mu_2 + \cdots +\mu_i \label{ehd equation}
\end{eqnarray} 
or inequality~(\ref{eh- equation}) fails.

First note that if $j=0$, then $S_n(\nu) = S_n(\nM)$, so 
inequality~(\ref{ehd equation}) cannot hold for $i=n$ since $\nM \leq \mu$.

Now suppose that there exists an $i\leq n-2$  such that 
inequality~(\ref{ehd equation}) holds.  Choose the smallest such $i$.
Then $\nu_i > \mu_i$, and both are odd integers, so $\nu_i -2 \geq \mu_i$.   
Suppose $\nu_i$ is in the $k^{\textrm{th}}$ batch of $\nu$.

 We consider the $(i+1)^{\textrm{th}}$ to $\Nk^{\textrm{th}}$ entries of 
$\nu$ and $\mu$.  Since $\nu$ is $M$-dominant and $M$-$2$-minuscule, 
$\nu_{i+1}, \ldots, \nu_{\Nk} \in \{\nu_i, \nu_i -2  \}$.  
Thus 
\begin{equation*}
\nu_{i+1}+ \cdots + \nu_{\Nk} \geq  (\Nk-i)(\nu_i-2).
\end{equation*} 
Also, since $\mu$ is $G$-dominant and $\mu_i \leq \nu_i -2$, it follows that 
$\mu_{i+1} + \cdots + \mu_{\Nk}  \leq  (\Nk-i)(\nu_i -2)$.  
Thus 
\begin{eqnarray*}
\mu_{i+1} + \cdots + \mu_{\Nk} &\leq & \nu_{i+1}+ \cdots + \nu_{\Nk}.
\end{eqnarray*}  
Combining this with inequality~(\ref{ehd equation}) yields 
\begin{eqnarray}
\mu_1 + \cdots + \mu_{\Nk} &<& \nu_1 + \cdots + \nu_{\Nk}. \label{ehF equation}
\end{eqnarray}

If $j \geq 2$ and $\nu_i$ is in the final batch, then we have 
$\mu_i < \nu_i = \pm 1$, so $\mu_i \leq -1$.  Since $\mu$ is $G$-dominant, 
it follows that $i=n$.  This contradicts our assumption that $i \leq n-2$.  
If $j \geq 2$ and $\nu_i$ is not in the final batch,
then $\sigma(k) \leq \sigma(r) \leq n-2$ and inequality~(\ref{ehF equation})
 contradicts $\nM\leq\mu$.

If $j=0$ and $\sigma(k)\not= n-1$, then inequality~(\ref{ehF equation}) 
contradicts 
$\nM\leq\mu$.
If $j=0$ and $\sigma(k) = n-1$, then we have $n_r=1$. 
As in the proof of Theorem~\ref{MTD}, we observe that 
inequality~(\ref{eh- equation}) holds, 
which when added to inequality~(\ref{ehn equation}) yields 
$S_{n-1}(\nu) \leq S_{n-1}(\mu)$,
contradicting inequality~(\ref{ehF equation}).

Now suppose that inequality~(\ref{eh- equation}) fails.  Thus we have
\begin{eqnarray}
S_{n-1}(\nu)-\nu_n &>& S_{n-1}(\mu)-\mu_n.  \label{ehi equation}
\end{eqnarray}
By what we have already shown, we have  
\begin{eqnarray}
S_{n-2}(\nu) &\leq & S_{n-2}(\mu). \label{ehC equation}
\end{eqnarray}
We claim that
\begin{eqnarray}
\nu_{n-1}-\nu_n &>& \mu_{n-1}-\mu_n \geq 0.  \label{ehj equation}
\end{eqnarray}
The first inequality follows from inequalities (\ref{ehi equation}) and
(\ref{ehC equation}), and the second holds since $\mu$ is $G$-dominant.
We have shown that inequality~(\ref{eh- equation}) holds if $j=0$ and $n_r=1$
 so we may assume that this is not the case.  Thus $\nu_{n-1}$ and $\nu_n$ 
are in the same batch of $\nu$, so, since $\nu$ is $M$-minuscule,
$\nu_{n-1}-\nu_n\in \{ 0, 2 \}$.
Combining this with  inequality~(\ref{ehj equation}) 
gives  
\begin{eqnarray}
\nu_{n-1}- \nu_n =2  & \textrm{and} & \mu_{n-1} - \mu_n =0. \label{ehD equation}
\end{eqnarray}
Substituting these values into 
inequality  (\ref{ehi equation}) and  combining the result with 
inequality~(\ref{ehC equation})  gives
\begin{eqnarray*}
S_{n-2}(\nu)+2>S_{n-2}(\mu)\geq S_{n-2}(\nu).  
\end{eqnarray*}
Both are integers of the same parity, so $S_{n-2}(\mu) =S_{n-2}(\nu)$.
This contradicts condition~(\ref{eh** equation}) since it follows from 
the equalities in (\ref{ehD equation}) (bearing in mind that all entries of
$\mu$ and $\nu$ are odd) that $S_{n}(\mu)$ is not congruent modulo four to 
$S_{n-2}(\mu)$, and $S_{n}(\nu)$ is congruent modulo four to  $S_{n-2}(\nu)$. 
 
Finally, suppose that inequality~(\ref{ehd equation}) holds for $i=n$.  We 
have already handled the $j=0$ case so we may assume that $j \geq 2$.  
Therefore $\nu_{n-1}$ and $\nu_n$ will be in the 
$r+1^{\textrm{st}}$ batch, so since $\nu$ is 
$M$-$2$-minuscule, we have $\nu_{n-1}+\nu_n \in \{ 0,2 \}$.  Also, since 
$\mu$ is $G$-dominant, we have $\mu_{n-1}+ \mu_n \geq 0$.  
 We have shown that $S_{n-2}(\nu) \leq S_{n-2}(\mu)$. 
It follows that  $S_n(\nu) \leq S_n(\mu)$ unless  
$S_{n-2}(\nu) = S_{n-2}(\mu)$, 
$\nu_{n-1}+\nu_n = 2$, and $\mu_{n-1}+ \mu_n = 0$.  In this case, 
$S_n(\nu) = S_n(\mu) +2$ which contradicts condition~(\ref{eh** equation}).
Thus $\nu\leq\mu$. 
\end{proof}

\end{document}